\documentclass[12pt,a4paper,twoside,final,notitlepage, leqno]{article}
\usepackage[english]{babel}
\usepackage[T1]{fontenc}  
\usepackage{graphicx}
\usepackage{amsmath,amsthm,epsfig,amsfonts,bbm}  

\newcommand{\pequationdeb}{$$ \left\{ \begin{minipage}[c]{130mm}}
\newcommand{\pequationfin}{\end{minipage}
                           \right. $$}

\def \smb {{\scriptstyle \bullet }}

\newcommand{\beq}     {\begin{equation}}
\newcommand{\enq}     {\end{equation}}
\newcommand{\be}    {\begin{enumerate}}
\newcommand{\ee}    {\end{enumerate}}

\newcommand{\Bb}

\newcommand{\R}{\mathbb{R}}

\newcommand{\sth}{\smash{^{^{\displaystyle \!\!\! \star}}}}

\def\resume{\if@twocolumn
\section*{R\'esum\'e}
\else \small
\quotation{\bf \it R\'esum\'e \rule[1mm]{1.5mm}{0.2mm}\vspace{0pt}}
\fi}
\def\endresume{\if@twocolumn\else\endquotation\fi}
\def\abstract{\if@twocolumn
\noindent\section*{{\bf Abstract}}
\else \small
\quotation{\noindent \bf {Abstract.} \rule[1mm]{1.5mm}{0.2mm}\vspace{0pt}}
\fi}
\def\endabstract{\if@twocolumn\else\endquotation\fi}

\def\section*#1{}
\begin{document}  

\title{\bf \Large      Petrov-Galerkin Finite Volumes   ~\\~  }

\author { { \large  Fran\c{c}ois Dubois~$^{a}$}   \\ ~\\
{\it  $^a$  Applications Scientifiques du Calcul Intensif, 
b\^at. 506, BP 167}, \\
{\it  F-91~403 Orsay  Cedex, Union Europ\'eenne.} \\ dubois@asci.fr }

\date { \large    26 april  2002~\protect\footnote{~  
{\rm  \small $\,\,$ Presented at the Third ``FVCA'' Conference, 
  Porquerolles (France),  24-28 June 2002.  Published  p.~203-210   in   
the Proceedings  {\it Finite Volumes for Complex Applications III; 
Problems \& Perspectives}, 
  Rapha\`ele Herbin  and  Dietmar Kr\"oner (Eds),   Hermes Prenton Press, May  2002.
Edition 7 November 2010.  }  }     }

\maketitle

%%%%%%%%%%%%%%%%%%%%%%%%%%%%%%%%%%%%%%%%%%%%%%%%%%%%%%%%%%%%%%%%%%%%%%%%%%%%%%%%%%%%%%%

\bigskip 

\noindent {\it \bf  Abstract } 

\noindent  
For an elliptic problem with two space dimensions, we propose to formulate the finite
volume method with the help of Petrov-Galerkin mixed finite elementsthat are based on
the building of a dual Raviart-Thomas basis.  

\smallskip \noindent {\it \bf   R\'esum\'e } 

\noindent   Pour un probl\`eme elliptique bidimensionnel, nous proposons de
formuler la m\'ethode des volumes   finis avec des \'el\'ements finis mixtes de
Petrov-Galerkin qui s'appuient sur la construction d'une base duale de
Raviart-Thomas.  

\smallskip \noindent {\it \bf   Keywords } 
  Finite volumes, mixed finite
elements, Petrov-Galerkin variational formulation, inf-sup condition, Poisson
equation.

%%%%%%%%%%%%%%%%%%%%%%%%%%%%%%%%%%%%%%%%%%%%%%%%%%%%%%%%%%%%%%%%%%%%%%%%%%%%%%%%%%%%%% 
 \bigskip    \bigskip  \noindent {\bf   \large  1. \quad   Introduction}  

\noindent $\bullet \quad$ 
 Let $\, \Omega \,$ be a bidimensional bounded convex domain
in $\, \R^2 \,$ with a polygonal boundary $\,\, \partial \Omega. \,$ We consider   the
homogeneous Dirichlet problem for the Laplace operator in the domain $\, \Omega \,$~:
\begin{equation}  \left\{ \begin{array}{ll}
- \Delta u \,\,=\,\, f \,  \qquad{\rm in}\,\,    \Omega  \,  \\
\,\,\,\,  u \,\,\,\,\,\,
\,\,=\,\, 0 \,  \qquad \,\, {\rm on \,\, the \,\, boundary  \,\, \partial
\Omega  \,\, of} \,\, \Omega  . \, 
\end{array} \right.\label{(1.1)} \end{equation}
We suppose that the datum $\, f \,$ belongs to the space $\, L^2(\Omega) .\,$ 
We introduce the momentum  $\,\, p \,\,$  defined by 
\begin{equation}
p \,\, = \,\, \nabla u \,.\,
\label{(1.5)} \end{equation}
Taking  the divergence of both terms arising in equation (\ref{(1.5)}), taking into
account the relation (\ref{(1.1)}), we observe that 
the divergence of momentum $\,p \,$ belongs to the space $\, L^2(\Omega) .\,$ For this
reason, we introduce the vectorial Sobolev space $\, H({\rm div},\, \Omega)  \,=\,
 \bigl\{ \,q \in  L^2(\Omega) \times L^2(\Omega) \,
,\,\, {\rm div} \, q \in  L^2(\Omega) \, \bigr\} .\,\, $
The variational formulation of the  problem (\ref{(1.1)}) with the help of the pair $\,
\xi \,=\, (u ,\, p) \,$ is obtained by testing the definition  (\ref{(1.5)}) against a vector
valued function $ \, q \,$ and integrating by parts. With the help of the boundary
condition, it comes : $ \,\, 
(p,\,q) \,+\, (u,\, {\rm div} \, q) \,=\, 0 ,\,\, \forall \, q \in 
 H({\rm div},\, \Omega) . \, \, $
Independently, the relations  (\ref{(1.1)}) and  (\ref{(1.5)}) are integrated on 
the domain $\, \Omega \,$ after multiplying by a scalar valued function $\, v \, 
\in  L^2(\Omega) \,.\,$
We obtain : $ \,\, 
({\rm div} \, p ,\,v ) \,+\, (f,\, v) \,=\, 0 ,\,\,  \forall \, v \in 
L^2(\Omega) .\,  \, $
The ``mixed'' variational formulation is obtained by introducing the product space
$\, V \,$ defined as  $ \,\, 
V \,=\, L^2(\Omega) \, \times \,  H({\rm div},\, \Omega) ,\,\, $ $ \,\,
\parallel (u,\,p)  \parallel_{V}^2 \,\,\equiv \, \parallel u
\parallel_{0}^2 \,\,+\,\, \parallel p \parallel_{0}^2 \,\,+\,\, \parallel {\rm
div}\,p  \parallel_{0}^2  ,\,\, $
the following bilinear form $\, \gamma(\smb,\, \smb) \,$ defined on $\, V \, \times \,
V \,: $ 
\begin{equation}
\gamma \bigl( (u,\,p),\, (v,\,q) \, \bigr) \,\,= \,\, (p,\,q) \,+\, (u,\, {\rm div}
\, q) \,+\,  ({\rm div} \, p ,\,v ) \,
 \label{(1.11)}  \end{equation} 
and the linear form $\, \sigma(\smb) \,$ defined on $\, V   \, $ according to : 
$ \,\, < \sigma ,\, \zeta > \,=\, -(f,\,v) ,\,\,, \zeta = (v,\, q) \, \in V . \,$
Then the Dirichlet problem (\ref{(1.1)}) takes the form : 
\begin{equation}  \left\{ \begin{array}{ll} 
\xi \, \in V  \\
\gamma(\xi,\, \zeta) \,\,=\,\, < \sigma ,\, \zeta > \,\,,\qquad \forall \, \zeta \,
\in V \,.\,
\end{array} \right. \label{(1.13)} \end{equation}
Due to classical inf-sup conditions introduced by Babu\v{s}ka \cite{Ba71}, 
the problem (\ref{(1.13)}) admits a unique solution $\, \xi \, \in V \,.\,$

%%%%%%%%%%%%%%%%%%%%%%%%%%%%%%%%%%%%%%%%%%%%%%%%%%%%%%%%%%%%%%%%%%%%%%%%%%%%%%%%%%%%%% 
\bigskip  \bigskip   \noindent {\bf   \large  2. \quad    Mixed Finite Elements}

\noindent $\bullet \quad$  
We introduce a mesh $\, {\cal T} \,$ that is   a bidimensional cellular complex
  composed in our case by triangular
elements $\, K \, $ $\, (K \, \in \,  {\cal E}_{\cal T} )  ,\,$  straight edges 
$\, a \, $ $\, (a \, \in \, {\cal A}_{\cal T}) \,$ and ponctual nodes $\, S \,$ 
$\, (S \, \in \, {\cal S}_{\cal T}) .\,$ We conside also classical finite
dimensional spaces $\, L_{\cal T}^2(\Omega) \,$ and
$\,  H_{\cal T}({\rm div},\, \Omega) \,$  that approximate respectively the spaces
$\, L^2(\Omega) \,$ and $\,  H({\rm div},\, \Omega) . \,$ A scalar valued function $\,
v \in  L_{\cal T}^2(\Omega) \,$ is constant in each triangle $\, K \,$ of the mesh : 
$ \,\, 
L_{\cal T}^2(\Omega) \,=\, \bigl\{ v : \Omega \longrightarrow \R ,\, \forall \, K 
\in {\cal E}_{\cal T} ,\, \exists \, v_K \, \in \R ,\, \forall \, x \in K, \, v(x)
\,=\, v_K \,\bigr\} \,.\,$
A vector valued function  function $\, q \in  H_{\cal T}({\rm div},\, \Omega) \,$ is
a linear combination of Raviart-Thomas \cite{RT77} basis functions $\, \varphi_a \,$ of
lower degree, defined for each edge $\, a \, \in {\cal A}_{\cal T} \,$ as follows. 

\noindent $\bullet \quad$  
Let $\, a \, \in {\cal A}_{\cal T} \,$ be an internal edge of the mesh, denote by
$\, S \,$ and $\, N \,$ the two vertices that compose its boundary (see Figure 1) :
$ \,\, \partial a \,=\, \{ \, S ,\, N \,\} \,$
and by $\, K \,$ and $\,  L \,$ the two elements that compose its co-boundary $\,
\partial^c a \,\equiv\, \{ \, K ,\, L \,\} \, $
in such a way that the normal direction $\, n \,$ is oriented from $\, K \, $ towards
$\, L \,$ and that the pair of vectors $\, (n,\, \,   {\overrightarrow {\rm SN}}) \,$
 is direct, as
shown on Figure 1. We denote by $\, W\,$ (respectively by $\, E )\,$ the third vertex
of the triangle $\, K \,$ (respectively of the triangle $\, L ) \,$ : 
$ \,\, K \,=\, (S,\, N ,\, W ), \,\, L \,=\, (N,\, S ,\, E ).\, $
% 

%%%%%%%%%%%%%%%%%%%%%%%%%%%%%%%%%     Figure 1        %%%%%%%%%%%%%%%%%%%%%%%%%%%%%% 
%
\bigskip 
\centerline { \includegraphics[width=5cm]{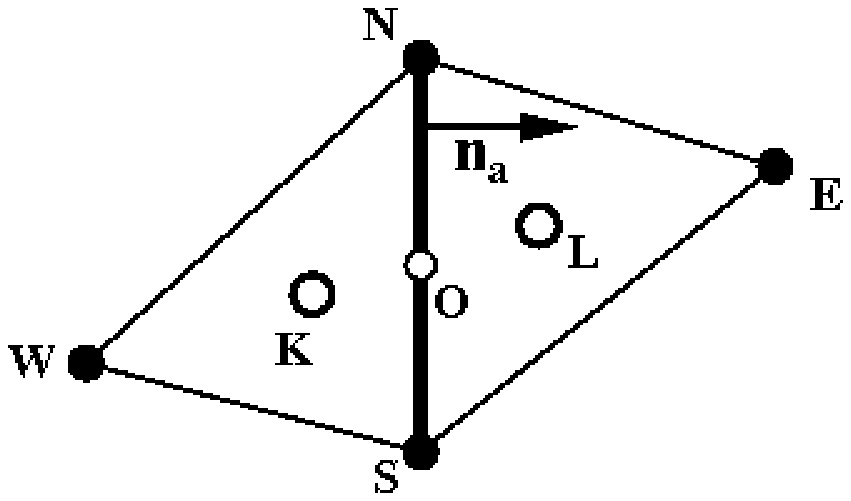} } 

\smallskip \centerline { {\bf Figure 1.} Co-boundary $ \, (K,\,L) \,$ of the edge $\,
a=(S,\,N).\,$  }
\smallskip

%%%%%%%%%%%%%%%%%%%%%%%%%%%%%%%%%%%%%%%%%%%%%%%%%%%%%%%%%%%%%%%%%%%%%%%%%%%%%%%%%%%%
%
\noindent
The vector valued   Raviart-Thomas \cite{RT77} basis  function $\,\varphi_a \,$ is defined
by the relations  $ \,\,  \varphi_a (x) \,\,= \,\,   \frac{1}{2 \mid K \mid } \, (x-W) 
\, {\rm when} \,  x \, \in K ,\, $  $   \varphi_a (x) \,\,= \,\, 
-\frac{1}{2 \mid L \mid } \, (x-E) 
\, {\rm when} \,  x \, \in L \, \, $ and $   \varphi_a (x) \,\,= \,\, 0 \, $  elsewhere.
When the edge $\, a \,$ is on the boundary $\, \partial \Omega $, we suppose that
the normal $ n $ points towards the exterior of the domain, so the element $\, L
\,$ is absent. We have in all cases  the $\,  H({\rm div},\, \Omega) \,$ conformity : 
$ \,\, \varphi_a  \in  H({\rm div},\, \Omega)  \,\,$
and the degrees of freedom are the fluxes of vector field $\, \varphi_a \,$ for all
the edges of the mesh (see \cite{RT77}) : $ \,\, 
\int_b \, \varphi_a  \, \smb \, n \,\, {\rm d}\gamma \,=\, \delta_{a,\,b} ,
\,\, \forall \, a ,\, b \, \in {\cal A}_{\cal T} .\, \,$
A vector valued function   $\, q \in  H_{\cal T}({\rm div},\, \Omega) \,$ is
a linear combination of the basis functions $\, \varphi_a \,: \,$ $ \,\, 
q \,=\, \sum_{a \, \in {\cal A}_{\cal T}} \, q_a \, \varphi_a \,    \in
\,   H_{\cal T}({\rm div},\, \Omega) \,=\, < \varphi_b \,,\,\, b \, \in
{\cal A}_{\cal T} > .\,\,$

 \noindent $\bullet \quad$ 
The mixed finite element method consist in choosing as discrete linear space the
following product : $ \,\, 
V_{\cal T} \,=\, L_{\cal T}^2(\Omega) \,\times \,  H_{\cal T}({\rm div},\,
\Omega) \,\,$
and to replace the letter $\, V \,$ by $\, V_{\cal T}  \,$ inside the variational
formulation   (\ref{(1.13)}) :  $ \,\, 
\xi_{\cal T}  \, \in V_{\cal T}   ,\,\, $ $ \,\, 
\gamma(\xi_{\cal T} ,\, \zeta) \,=\, < \sigma ,\, \zeta > ,\,\, \forall \,
\zeta \, \in V_{\cal T}  .\, \,$ 
In other terms 
\begin{equation}  \left\{ \begin{array}{ll} 
  u_{\cal T} \, \in  L_{\cal T}^2(\Omega)
\,\,,\quad  p_{\cal T} \, \in   H_{\cal T}({\rm div},\, \Omega) \,  \\ 
(p_{\cal T} \,,\, q ) \,+\,
(u_{\cal T} \,,\, {\rm div} \, q  ) \,\,=\,\, 0 \,,\qquad \forall \,
q \, \in   H_{\cal T}({\rm div},\, \Omega) \, \\ 
({\rm div} \, p_{\cal T} \,,\, v )
\,+\, (f\,,\, v ) \,\, \,\,\,\,=\,\, 0 \,,\qquad \forall \, v  \,
\in  L_{\cal T}^2(\Omega) \,. \,
\end{array} \right. \label{(1.24)} \end{equation} 
The numerical analysis of the relations between the continuous problem  (\ref{(1.13)})
and the discrete problem (\ref{(1.24)}) as the mesh $\,{\cal T}\,$  is more and more
 refined is
classical \cite{RT77}. The above method is popular in the context of petroleum and nuclear
industries   but suffers from the fact that the associated linear system is
quite difficult to solve from a practical point of view. The introduction of
supplementary Lagrange multipliers by Brezzi, Douglas and Marini \cite{BDM85} allows a
simplification of these algebraic aspects, and  their interpretation by Croisille in
the context of box schemes \cite{Cr2k} gives a good mathematical foundation of a popular
numerical method.

%%%%%%%%%%%%%%%%%%%%%%%%%%%%%%%%%%%%%%%%%%%%%%%%%%%%%%%%%%%%%%%%%%%%%%%%%%%%%%%%%%%%%% 
\bigskip  \bigskip  \noindent {\bf   \large  3. \quad     Finite Volumes}

 \noindent $\bullet \quad$ 
From a theoretical and practical  point of view, the resolution of the linear system 
(\ref{(1.24)}) can be conducted as follows. We introduce the mass-matrix 
$\,\, M_{a,\,b} \, = \, (\varphi_a ,\, \varphi_b) ,\, \,  a,\, b \, \in {\cal
A}_{\cal T} \, \,$ associated with the Raviart-Thomas vector valued functions.
Then the first equation of (\ref{(1.24)}) determines the momentum  $ \,\,
p_{\cal T} \, = \,  \sum_{a \, \in   {\cal A}_{\cal T}} $ $ \,p_{{\cal
T}\!\!,a} \, \varphi_a \,\, $
as a function of the mean values $\,\, u_{{\cal T}\!\!,K} \, $ for $\, K \, \in 
{\cal E}_{\cal T} \,:\,$ 
\begin{equation} 
p_{{\cal T}\!\!,a} \,\,=\,\,- \sum_{b \, \in  {\cal A}_{\cal T}}\,  \bigl( M^{-1}
\bigr)_{a,\,b}  \,\,  \sum_{ K \, \in  {\cal E}_{\cal T}} \, u_{{\cal T}\!\!,K}
\,\, \int_K \, {\rm div}\, \varphi_b \, {\rm d}x \,\,.  
\label{(1.27)} \end{equation}  
The representation (\ref{(1.27)}) suffers at our opinion form a major defect : due to the
fact that the matrix $\,\,  M^{-1} \,\,$ is full, the discrete gradient $
\,p_{\cal T} \,$ is a {\bf global} function of the mean values $\,  u_{{\cal
T}\!\!,K} \,$ and this property contradicts the mathematical  foundations of the
derivation operator to be  {\bf linear} and {\bf local}. An a posteriori
correction of this defect has been proposed by Baranger, Ma\^{\i}tre and Oudin 
\cite{BMO96} : with an appropriate  numerical integration of the mass matrix $\, M ,\,$ 
it is possible to lump it and the discrete gradient in the direction $\, n \,$ of
the edge $\, a \,$ is represented by a formula of the type~: 
\begin{equation} 
p_{{\cal T}\!\!,a} \,\,=\,\,\frac {u_{{\cal T}\!\!,L} - u_{{\cal T}\!\!,K}}{h_a} \,
\label{(1.28)} \end{equation}  
with the notations of Figure 1. The substitution of the relation (\ref{(1.28)}) inside the
second equation of the formulation (\ref{(1.24)}) conducts to a variant of the so-called
finite volume method. In an analogous manner, the family of finite volume schemes
proposed by Herbin \cite{He95} suppose {\it a priori} that the discrete gradient in the
normal direction admits a representation of the form (\ref{(1.28)}). Nevertheless, the
intuition is not correctly satisfied by a scheme such that (\ref{(1.28)}). The finite
difference $\,\, \smash{\frac{u_{{\cal T}\!\!,L} - u_{{\cal T}\!\!,K}}{h_a}} \,\,$
wish to be a  a good approximation of the gradient $\,\, p_{\cal T} \,=\, \nabla
u_{\cal T} \,\,$ in the direction $ {\overrightarrow {\rm KL}}   $ whereas
 the coefficient $\,\,
p_{{\cal T}\!\!,a} \,\,$ is an approximation of $\,\,  \int_{\rm \displaystyle a}
\nabla  u_{\cal T} \smb n\,\, {\rm d}\tau \,\,$ in the {\bf normal} direction (see
again the Figure 1). When the mesh $\, {\cal T} \,$ is composed by general triangles,
this approximation is not completely satisfactory and contains  a real limitation of
these variants of the finite volume method at our opinion. 

 \noindent $\bullet \quad$ 
In fact, the finite volume method for the approximation of the diffusion operator
 has been first proposed from empirical considerations. Following {\it e.g.}
Noh \cite{No64} and Patankar \cite{Pa80}, the idea is to represent the normal interface 
gradient  $\,\,  \int_{\rm \displaystyle a}   \nabla  u_{\cal T} \smb n\,\, {\rm d}  \tau
\,\,$ as a function of {\bf neighbouring} values. Given an edge $\, a,\,$ a  vicinity
$\, {\cal V}(a) \,$ is first determined in order to represent the normal gradient
$\,\,p_{{\cal T}\!\!,a}  =  \int_{\rm \displaystyle a}   \nabla  u_{\cal T} \smb
n\,\, {\rm d} \tau \,\,$  with a formula of the type 
\begin{equation}
 \int_{\rm \displaystyle a}   \nabla  u_{\cal T} \smb n\,\, {\rm d} \tau \,\,=\,\,
 \sum_{K  \in {\cal V}(a)} \, g_{a,K} \,\, u_{{\cal T}\!\!,K}  \,. \,
\label{(1.29)} \end{equation}  
Then the conservation equation  $ \,\, 
{\rm div} \, p \,+\, f \,=\, 0 \,\,$
is integrated inside each cell $\, K  \in {\cal E}_{\cal T} \,$ is order to
determine an equation for the mean values $\,  u_{{\cal T}\!\!,K} \,\,$ for all $\,
K  \in {\cal E}_{\cal T} .\,$ The difficulties of such approches have been presented
by Kershaw \cite{Ke81} and a variant of such  scheme has been   first analysed by
Coudi\`ere, Vila and Villedieu \cite{CVV99}. The key remark that we have done with
F.~Arnoux (see \cite{Du89}), also observed by Faille, Gallou\"et and Herbin 
\cite{FGH91}  is that the representation (\ref{(1.29)}) must be 
{\bf exact  for linear} functions $\, u_{\cal
T} \,.\,$ We took this remark as a starting point for our tridimensional 
finite volume scheme proposed in \cite{Du92}. It is also an essential hypothesis for the
result proposed by Coudi\`ere, Vila and Villedieu.

%%%%%%%%%%%%%%%%%%%%%%%%%%%%%%%%%%%%%%%%%%%%%%%%%%%%%%%%%%%%%%%%%%%%%%%%%%%%%%%%%%%%%% 
\bigskip  \bigskip \noindent {\bf   \large  4. \quad   
Finite volumes as mixed Petrov-Galerkin  finite elements } 

 \noindent $\bullet \quad$ 
In this contribution, we propose to  discretize the variational problem  (\ref{(1.13)})  
with the Petrov-Galerkin mixed finite element method, first introduced by Thomas and
Trujillo [TT99]. In the way  we have proposed  in \cite{Du2k}, the idea is to
construct a  discrete  functional  space 
$\,\,  H_{\cal T}\sth({\rm div},\, \Omega) \, \,$   generated by vectorial functions 
$\, \varphi^{\star}_{a} ,\,  a   \in {\cal A}_{\cal T} ,\,$ 
that are conforming in the space $\, H({\rm div},\, \Omega) \, :\,$ $ \,\,  
\varphi^{\star}_{a} \, \in \, H({\rm div},\, \Omega) \, $
and to represent exactly  the {\bf dual basis} of the family $\, \{ \varphi_{b},\, b 
\in {\cal A}_{\cal T} \}\,$  with the $\, L^2 \,$ scalar product~: $ \,\, 
( \, \varphi_{a} \,,\,\varphi_{b}^{\star} \, ) \,\,= \,\, \delta_{a,\,b} ,\,\, 
\forall \, a ,\, b \, \in {\cal A}_{\cal T} .  \,\, $ Then $ \,\, 
H_{\cal T}\sth({\rm div},\, \Omega) \,=\, <
\varphi_{b}^{\star}  \,,\,\, b \, \in {\cal A}_{\cal T} > .\,\,$
The mixed Petrov-Galerkin mixed finite element method consists just in
replacing  the space $\,\, H_{\cal T} ({\rm div},\, \Omega) \,\,$ by the dual space 
$\,\, H_{\cal T}\sth({\rm div},\, \Omega) \,\,$  for {\bf test functions} in the
first equation of discrete formulation (\ref{(1.24)}). We obtain by doing this the so-called
{\bf  Petrov-Galerkin finite volume} scheme : 
\begin{equation}  \left\{ \begin{array}{ll}   
u_{\cal T} \, \in  L_{\cal T}^2(\Omega) \,\,,\quad  p_{\cal T} \, \in 
  H_{\cal T}({\rm div},\, \Omega) \, \\ 
(p_{\cal T} \,,\, q ) \,+\,
(u_{\cal T} \,,\, {\rm div} \, q  ) \,\,=\,\, 0 \,,\qquad \forall \,
q \, \in    H_{\cal T}\sth({\rm div},\, \Omega)  \, \\
 ({\rm div} \, p_{\cal T} \,,\, v )
\,+\, (f\,,\, v ) \,\, \,\,\,\,=\,\, 0 \,,\qquad \forall \, v  \,
\in  L_{\cal T}^2(\Omega) \,. \,
\end{array} \right. \label{(1.34)} \end{equation}
We introduce a compact form of the previous mixed Petrov-Galerkin formulation with the
help of the product space $\, V_{\cal T}^{\displaystyle \star} \,$ defined by $ \,\, 
V_{\cal T}^{\displaystyle \star} \,=\, L_{\cal T}^2(\Omega) \,\times \, 
H_{\cal T}\sth ({\rm div},\, \Omega) . \,\,$
Then the discrete variational formulation (\ref{(1.34)}) admits the form : $ \,\,
\xi_{\cal T}  \, \in V_{\cal T}  , \,\, $ $ \,\, 
\gamma(\xi_{\cal T} ,\, \zeta) \,=\, < \sigma ,\, \zeta > , \, \forall \,
\zeta \, \in V_{\cal T}^{\displaystyle \star}  . \,\, $

%%%%%%%%%%%%%%%%%%%%%%%%%%%%%%%%%%%%%%%%%%%%%%%%%%%%%%%%%%%%%%%%%%%%%%%%%%%%%%%%%%%%%% 
\bigskip  \bigskip \noindent {\bf   \large  5. \quad     Stability analysis}

\noindent $\bullet \quad$ 
We suppose in the following that the mesh $\, \, {\cal T} \,\,$ is a bidimensional 
cellular complex composed by triangles as proposed in the previous sections. 
 Following the work of Ciarlet and Raviart \cite{CR72}, for
any element $\,\, K   \in \,  {\cal E}_{\cal T} \, \,$ we denote by $\,h_{\!K}\,$
the diameter of the triangle $\, K \,$ and by $\, \rho_{\!K}\,$ the diameter of the
inscripted ball inside $\, K .\,$ We suppose that the mesh $\, \, {\cal T} \,\,$
belongs   to a family $\, {\cal U}_{\theta} \,$ $\, (\theta > 0)\,$
of meshes  defined by the condition $ \,\, 
{\cal T} \,\in \, {\cal U}_{\theta} \,\,\, \Longleftrightarrow \,\,\, 
\forall \, K \, \in {\cal E}_{\cal T} \,,\,\, {{h_{\!K}}\over{\rho_{\!K}}} \,
\leq \, \theta . \,$ 
We suppose also that the dual space $\,\,  H_{\cal T}\sth({\rm div},\, \Omega) \,\,$
constructed by the previous conditions satisfies the following hypothesis. 

\noindent {\bf Hypothesis 1.} \quad {\bf Interpolation operator}
 $\,\, \,  H_{\cal T}({\rm div},\,
\Omega) \, \longrightarrow \, H_{\cal T}\sth ({\rm div},\,
\Omega) \,. \,$  
We suppose that the mesh $\, {\cal T} \,$ belongs to the family $\, {\cal
U}_{\theta} .\,$   Let $\,\, H_{\cal T}({\rm div},\, \Omega) \, \ni
\, q \,\, \longmapsto \,\, \Pi \, q \,\in \, H_{\cal T}\sth ({\rm div},\, \Omega)
\,\,$ be the mapping defined by the condition
\begin{equation}
\Pi \,\, \Bigl( \, \sum_{a \, \in {\cal A}_{\cal T}} \, q_a \,\, \varphi_a \,\Bigr)
\,=\, \sum_{a \, \in {\cal A}_{\cal T}} \, q_a \,\, \varphi^{\displaystyle 
\star}_a  ,\,\qquad  \, \sum_{a \, \in {\cal A}_{\cal T}} \, q_a \,\,
\varphi_a \,\,\, \in \,  H_{\cal T}({\rm div},\, \Omega) . \,\,
\end {equation}
We suppose that the dual basis $\,\, \varphi^{\displaystyle  \star}_a  \,\,$  is
constructed in such a way that there exists strictly positive constants $\, A,\, B,\,
D  ,\, E \,$ that only depends on the parameter $\, \theta \,$ such that we have the
following estimations : 
\begin{equation}  \left\{ \begin{array}{ll} 
A \, \parallel q \parallel_{0}^2 \,\,\,\, \leq \,\, ( \, q  \, ,\, \Pi \, q \, )
\,\,, \qquad \forall \, q \, \in  H_{\cal T}({\rm div},\, \Omega) \, \\
\parallel \Pi \, q \parallel_{0} \,\,\,\,\, \leq \,\, B \, \parallel  q
\parallel_{0} \,\,, \qquad \forall \, q \, \in  H_{\cal T}({\rm div},\, \Omega) \, \\ 
\parallel {\rm div} \,\Pi \, q \parallel_{0}  \,\,\,\, \leq \,\, D \,
\parallel {\rm div} \, q \parallel_{0}  \,\,, \qquad \qquad 
\forall \, q \, \in  H_{\cal T}({\rm div},\, \Omega) \, \\
( \, {\rm div}\, q \,,\,  {\rm div}\, \Pi \,q \,) \,\,\, \geq \,\, E \, 
\parallel {\rm div} \, q \parallel_{0}^2 \,\,, \qquad \forall \, q \, \in  H_{\cal
T}({\rm div},\, \Omega) \,. \, 
\end{array} \right. \label{(2.4)} \end{equation}

\noindent
{\bf Proposition 1.} \quad {\bf Technical lemma about lifting of scalar fields.} 
Let $\, \theta \,$ be a strictly positive parameter. We suppose that the dual
Raviart-Thomas basis satisfies the Hypothesis 1. Then there exists some strictly
positive constant $\, F \,$ that only depends on the parameter $\, \theta \,$ such that
for any mesh   $\, {\cal T} \,$  that belongs to the family $\, {\cal U}_{\theta} ,
\,$ and for any scalar field $\, u \,$ constant in each element $\, K  \in {\cal
E}_{\cal T} \,$ $(u \, \in \, L_{\cal T}^2(\Omega) ) ,\, $
there exists some vector field $\,\, q \, \in  \, H_{\cal T}\sth ({\rm div},\,
\Omega) \,\,$ such that $ \,\, 
\parallel q \parallel  _{H({\rm div},\, \Omega)} \,\leq \,F \, \parallel
u \parallel _{0} \, \, $ and $ \,\,  
(\, u \,,\, {\rm div} \, q \, ) \,\,  \geq \, \parallel u \parallel _{0}
^2 . \,\,  $ 

\noindent
{\bf Proposition 2} \quad {\bf Discrete stability. } 
Let $\, \theta \,$ be a strictly positive parameter. We suppose that the dual
Raviart-Thomas basis satisfies the Hypothesis 1. Then we have the following discrete
stability for the Petrov-Galerkin mixed formulation (\ref{(1.34)}) : $ \,\, 
\exists \, \beta   > \,0 \,,\,\,\, \forall \, {\cal T} \in \, {\cal U}_{\theta}
,\, \forall \, \xi \, \in V_{\cal T} \,\, $ such that $   \,\,\parallel \xi \parallel_{V}
  \,=\, 1 , \,$ $ \,\, \exists \, \eta \, \in  \, V_{\cal T}^{\displaystyle \star} ,\,\,
\parallel \zeta \parallel_{V}  \,\,\leq\,\, 1 \,,\, \, \gamma(\xi,\, \zeta) \,\,
\geq \,\, \beta ,\,\, $
with $\, \gamma (\smb,\, \smb) \,$ defined at the relation (\ref{(1.11)}) and 
$\, \beta \,$ chosen such that  $ \,\, 
\sqrt{ 1 - {{B + 2 D}\over{A}} \beta \,-\, \beta^2} \,\,\, \geq \,\,\, \Bigl( \, 
1 + F \, \bigl( 1 \,+\, \sqrt{{{B + 2 A }\over{A}}} \, \bigr) \, \Bigr) \,
\sqrt{\beta} . \, \, $

\smallskip \noindent
{\bf Theorem 1} \quad {\bf Optimal error estimate.} 
Let $\, \Omega \,$ be a two-dimensional open convex domain of $\, \R^2 \,$ with a
polygonal boundary,  $\, u \, \in H^2(\Omega) \,$ be the solution of the problem
(\ref{(1.1)})  considered under  variational formulation and $\, p = \nabla u \,$ be the
associated momentum. Let $\, \theta \,$ be a strictly positive parameter and  $\,
{\cal U}_{\theta} \,$  a family of meshes $\, {\cal T} \,$ that satisfy the
Hypothesis 1. Let $\, \xi\equiv (u_{\cal T},\,p_{\cal T}) \in V_{\cal T} \,$ be the 
solution of
the discrete problem (\ref{(1.34)}). Then there exists some constant $\, C > 0 \,$
that only depends on the parameter $\, \theta \,$ such that  $ \,\, 
\parallel u - u _{\cal T} \parallel _{0} \,+\, $ $ \parallel p - p _{\cal T}
\parallel  _{H({\rm div},\, \Omega)} \,    \leq  \, C \, \, h_{\cal T}
\,\,  \parallel f \parallel _{0} .\,\, $

%%%%%%%%%%%%%%%%%%%%%%%%%%%%%%%%%%%%%%%%%%%%%%%%%%%%%%%%%%%%%%%%%%%%%%%%%%%%%%%%%%%%%% 
\bigskip  \bigskip   \noindent {\bf   \large  6. \quad   
Towards a first Petrov-Galerkin finite volume scheme}

%%%%%%%%%%%%%%%%%%%%%%%%%%%%%%%%%     Figure 2        %%%%%%%%%%%%%%%%%%%%%%%%%%%%%% 
%
\bigskip 
\centerline { \includegraphics[width=6cm]{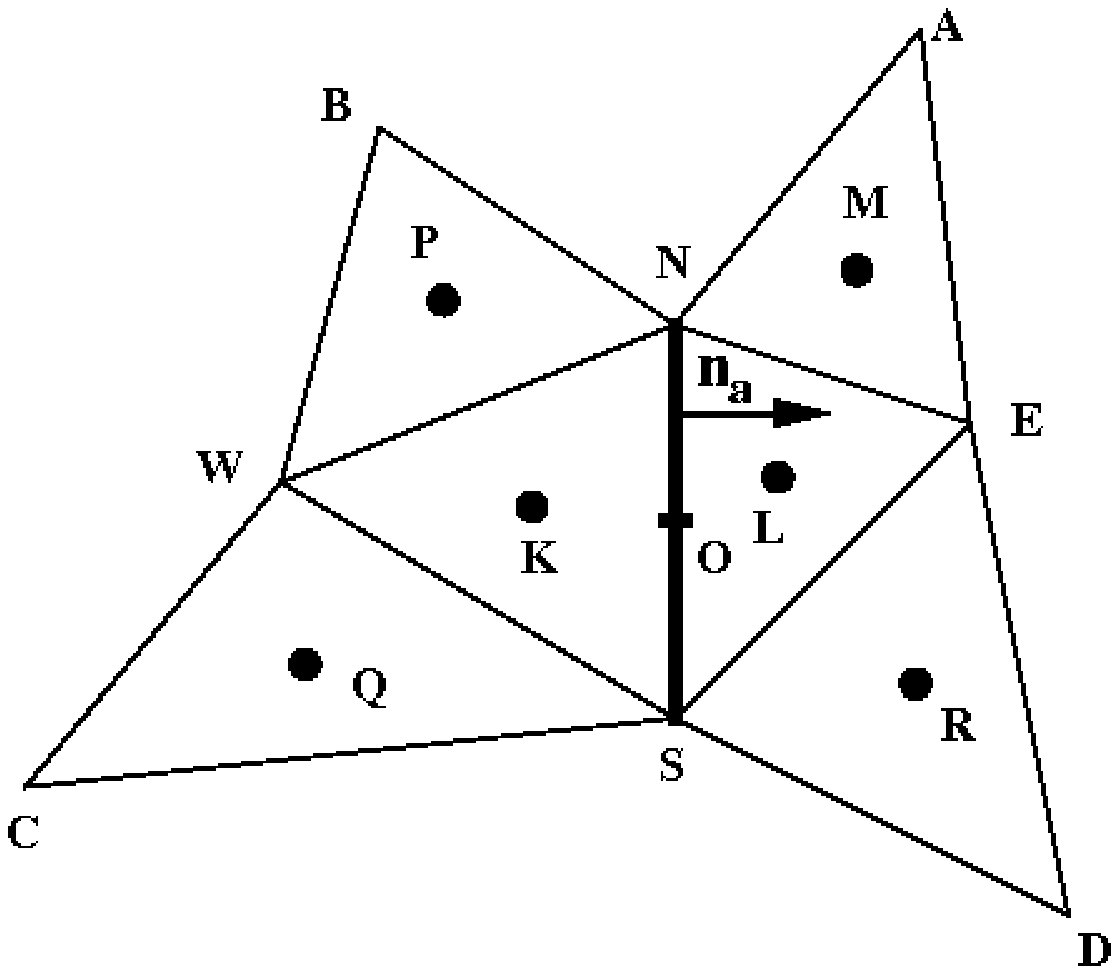} } 

\smallskip \centerline { {\bf Figure 2.}  Support $\, {\cal V}({\rm SN}) \,$ of the dual
Raviart-Thomas basis function $\, \varphi^{\displaystyle \star}_{\rm SN} .\,$ } 
\smallskip
%%%%%%%%%%%%%%%%%%%%%%%%%%%%%%%%%%%%%%%%%%%%%%%%%%%%%%%%%%%%%%%%%%%%%%%%%%%%%%%%%%%%
%

\noindent
{\bf Theorem 2} \quad We suppose that the internal edge $\, a \,$
links the two vertices $\, {\rm S }\,$ and $\, {\rm N} \,$ (see the Figure where  $\,
a = {\rm SN} , \,$ $\, {\rm O} \,$ is the middle of $\, {\rm SN} \,$ 
and  $\, n \,$ is the associated normal direction), if the
support of the dual Raviart-Thomas basis function $\, \varphi^{\displaystyle \star}_{\rm SN}
 \,$ is the  vicinity $\, {\cal V}(a) = \{  K,\,L,\, M,\, P,\, Q ,\, R \}
\,$ of the edge  $\, a \,$ composed by  six triangles presented on  Figure 2  
and if the divergence of the  dual Raviart-Thomas basis function is equal to a {\bf
constant field} in each triangle of $\, {\cal V}(a) \,$ $({\rm div} \, 
 \varphi^{\displaystyle \star}_{\rm SN}  \in   L_{\cal T}^2(\Omega)) ,\,$ then the 
{\bf five} mean flux values  $\, \eta \equiv
\int_{\rm SN} \,  \varphi^{\displaystyle \star}_{\rm SN}
\, \smb \, n \, {\rm d}\tau ,\,$ $\, \alpha
\equiv \int_{\rm EN} \, \varphi^{\displaystyle \star}_{\rm SN}
\, \smb \, n_{\rm EN} \, {\rm d}\tau
\,,\,$ $\, \beta \equiv  \int_{\rm WN} \, \varphi^{\displaystyle \star}_{\rm SN}
\, \smb \, n_{\rm
WN} \, {\rm d}\tau \,,\,$ $\, \gamma \equiv \int_{\rm WS} \, 
 \varphi^{\displaystyle \star}_{\rm SN} \, \smb
\, n_{\rm WS} \, {\rm d}\tau \,$ and $\, \delta \equiv \int_{\rm SE} \,
 \varphi^{\displaystyle \star}_{\rm SN} \, \smb \, n_{\rm SE} 
\, {\rm d}\tau \,$   satisfy the
following {\bf three scalar constraints} :  
\begin{equation}
\eta  \,   {\overrightarrow {\rm KL}}    \,+\, 
\alpha  \, {\overrightarrow {\rm LM}}    \,+\, 
\beta  \,  {\overrightarrow {\rm PK}}    \,+\, 
\gamma  \, {\overrightarrow {\rm QK}}    \,+\,
\delta  \, {\overrightarrow {\rm LR}}   \,\,=\,\, \mid {\rm {\overrightarrow
{\rm SN}}}  \mid \,  n \, 
\label{(2)} \end{equation} 
\begin{equation}
\alpha  \,  {\overrightarrow {\rm LM}}  \smb   {\overrightarrow {\rm WA}}  +  
\beta   \,  {\overrightarrow {\rm PK}}  \smb   {\overrightarrow {\rm EB}}  +  
\gamma  \,  {\overrightarrow {\rm QK}}  \smb   {\overrightarrow {\rm EC}}  + 
\delta  \,  {\overrightarrow {\rm LR}}  \smb   {\overrightarrow {\rm WD}}    
  =   -3   \mid {\rm {\overrightarrow {\rm SN}}} \mid n   \smb   \bigl( 
  {\overrightarrow {\rm OL}} +  {\overrightarrow {\rm OK}}   \bigr)  . \, 
\label{(3)} \end{equation} 

\noindent $\bullet \quad$ The finite volume approach  is then
obtained in the spirit of (\ref{(1.29)})
with a six point scheme for the mean gradient in the normal direction thanks to
the first equation of the mixed variational formulation (\ref{(1.34)}) : 
\begin{equation}
\int_{\rm SN}  \!\!\! \nabla  u_{\cal T} \smb n \, {\rm d}\tau =
\eta (u_{\rm L}\!-\!u_{\rm K}) + 
\alpha (u_{\rm M} \!-\! u_{\rm L}) + 
\beta (u_{\rm K} \!-\! u_{\rm P})  + 
\gamma (u_{\rm K} \!-\! u_{\rm Q}) +  
\delta (u_{\rm R} \!-\! u_{\rm L}) .\,
\label{(4)} \end{equation} 
We remark that the constraints (\ref{(2)}) express that the relation (\ref{(4)})
 is {\bf exact} if the field $\, u_{\cal T} \,$ is an affine function.   Taking into 
account the fact that we have five parameters for the definition of the finite 
volume scheme (relation   (\ref{(4)})) and only
three constraints (relations (\ref{(2)}) and (\ref{(3)})) for these parameters, 
the stability seems a reasonable goal, 
even if the problem  remains essentially open for general triangular meshes. 

% $\hfill$  26 april  2002.

%%%%%%%%%%%%%%%%%%%%%%%%%%%%%%%%%%%%%%%%%%%%%%%%%%%%%%%%%%%%%%%%%%%%%%%%%%%

\end{document}